\documentclass[11pt]{article}

\usepackage{amssymb}
\usepackage{amsmath}
\usepackage{amsfonts}
\usepackage{epsfig}
\usepackage{graphicx} 
\usepackage{multicol} 
\def\be{\begin{eqnarray}}
\def\ee{\end{eqnarray}}
\def\ben{\begin{eqnarray*}}
\def\een{\end{eqnarray*}}
\def\Ebox#1#2{%
\begin{center}
\parbox{#1\hsize}{\epsfxsize=\hsize \epsfbox{#2}}
\end{center}}
\def\ddt{\frac{d}{dt}}
\def\bfme{{\mbox{\protect\boldmath$e$}}} 
\def\bfms{{\mbox{\protect\boldmath$s$}}} 
\def\bfmy{{\mbox{\protect\boldmath$y$}}} 
\def\bfmB{{\mbox{\protect\boldmath$B$}}} 
\def\bfmN{{\mbox{\protect\boldmath$N$}}} 
\def\bfmS{{\mbox{\protect\boldmath$S$}}}
\def\bfmX{{\mbox{\protect\boldmath$X$}}}

\def\eqdef{\mathbin{:=}}

\def\sq{\hbox{\rlap{$\sqcap$}$\sqcup$}}\def\qed{\ifmmode\sq\else{\unskip\nobreak\hfil
\penalty50\hskip1em\null\nobreak\hfil\sq
\parfillskip=0pt\finalhyphendemerits=0\endgraf}\fi}

\def\proof{\paragraph{\sc Proof. }}

\newtheorem{theorem}{Theorem}[section]

\def\Theorem#1{Theorem~\ref{t:#1}}

\def\trace{\mathrm{ trace\,}} 

\def\transpose{{\hbox{\it\tiny T}}}
\def\barM{{\overline{M}}}
\def\barW{{\overline{W}}}
\def\barO{{\overline{O}}}

\def\Re{\field{R}}

\def\Co{\field{C}}
\def\eq#1/{(\ref{e:#1})}
\newcommand{\Th}[1]{Theorem~\ref{t:#1}}
\def\clL{{\cal L}} 
\def\lmin{\lambda_{\hbox{\rm\tiny min}}} 

\def\state{\mathsf{X}} 
\def\bx{\clB(\state)}
\def\La{{\clL_{\alpha}}}

\def\clQ{{\cal Q}}
\def\Qa{{\clQ_{\alpha}}}
\def\ha{h_\alpha}
\def\ho{h_0}
\def\la{\lambda_\alpha}
\def\lo{\lambda_0}

\def\LV{L_\infty^V}

\def\ddt{\frac{d}{dt}}

\def\Ebox#1#2{%
\begin{center}    
 \parbox{#1\hsize}{\epsfxsize=\hsize \epsfbox{#2}}
\end{center}}

\def\tlabel#1{\label{t:#1}}
\def\slabel#1{\label{s:#1}}
\def\flabel#1{\label{f:#1}}
\def\elabel#1{\label{e:#1}}

\def\sq{\hbox{\rlap{$\sqcap$}$\sqcup$}}
\def\qed{\ifmmode\sq\else{\unskip\nobreak\hfil
\penalty50\hskip1em\null\nobreak\hfil\sq
\parfillskip=0pt\finalhyphendemerits=0\endgraf}\fi}

\newsavebox{\junk}
\savebox{\junk}[1.6mm]{\hbox{$|\!|\!|$}}
\def\lll{{\usebox{\junk}}}

\def\limsup{\mathop{\rm lim\ sup}}

\def\state{{\sf X}}

\newcommand{\field}[1]{\mathbb{#1}}

\def\Re{\field{R}}

\def\Co{\field{C}}

\def\ind{\field{I}}

\newcommand{\RL}{{\mathbb R}}

\def\bfmx{{\mbox{\protect\boldmath$x$}}}

\def\bfmM{{\mbox{\protect\boldmath$M$}}}
\def\bfmW{{\mbox{\protect\boldmath$W$}}}

\def\bfPhi{\mbox{\protect\boldmath$\Phi$}}

\def\bfphi{\mbox{\boldmath$\phi$}}

\def\hatheta{{\hat\theta}}

\def\til={{\widetilde =}}

\def\tiltheta{{\tilde \theta}}

\def\clB{{\cal B}}

\def\clL{{\cal L}}
\def\clM{{\cal M}}

\def\clQ{{\cal Q}}

\def\eqdef{\mathbin{:=}}

\def\Prob{{\sf P}}

\def\Expect{{\sf E}}

\def\epsy{\varepsilon}

\def\varble{\,\cdot\,}

\def\Theorem#1{Theorem~\ref{t:#1}}

\def\Section#1{Section~\ref{s:#1}}

\def\eq#1/{(\ref{e:#1})}

\newcommand{\beqn}[1]{\notes{#1}%
\begin{eqnarray} \elabel{#1}}

\newcommand{\eeqn}{\end{eqnarray} } 

\newcommand{\beq}[1]{\notes{#1}%
\begin{equation}\elabel{#1}}

\newcommand{\eeq}{\end{equation}} 

\def\bdes{\begin{description}}
\def\edes{\end{description}}

\newcommand{\oo}{\overline}

\def\barf{{\oo {f}}}

\def\proof{\paragraph{\sc Proof. }}

\def\notes#1{}

\def\hatheta{\widehat{\theta}}
\def\tiltheta{\widetilde{\theta}}

\begin{document}
\title{The ODE Method \protect\newline 
  and Spectral Theory of Markov Operators}

\author
{
	J. Huang\thanks{Department of Electrical and Computer
                Engineering and the Coordinated Sciences Laboratory,
                University of Illinois at Urbana-Champaign,
                Urbana, IL 61801, U.S.A. Email: 
		{\tt jhuang@control.csl.uiuc.edu}.
		}
\and
        I. Kontoyiannis\thanks{Division
                of Applied Mathematics and Department
                of Computer Science, Brown University,
                Box F, 182 George St., Providence, RI 02912, USA.
                Email: {\tt yiannis@dam.brown.edu}
                Web: {\tt www.dam.brown.edu/people/yiannis/}.
                              }
\and
        S.P. Meyn\thanks{
		Department of Electrical and Computer
                Engineering and the Coordinated Sciences Laboratory,
                University of Illinois at Urbana-Champaign,
                Urbana, IL 61801, U.S.A. Email: {\tt s-meyn@uiuc.edu}.
                        }
}

\maketitle
 
\begin{abstract}
We give a development of the ODE method 
for the analysis of recursive algorithms
described by a stochastic recursion.
With variability modelled via an underlying 
Markov process, and under general assumptions,
the following results are obtained:
\begin{description}
\item[(i)]
Stability of an associated ODE implies that 
the stochastic recursion is stable in a strong 
sense when a gain parameter is small.
\item[(ii)]
The range of gain-values is quantified through a spectral analysis of
an associated linear operator, providing a non-local theory.
\item[(iii)]
A second-order analysis shows precisely how variability leads to
sensitivity of the algorithm with respect to the gain parameter.
\end{description}
All results are obtained within the
natural operator-theoretic framework
of geometrically ergodic Markov processes.
\end{abstract}

\section{Introduction}
\label{Introduction}

Stochastic approximation algorithms and their variants are commonly
found in control,  communication and related fields.  Popularity has
grown due to increased computing power, and the interest in various
`machine learning' algorithms \cite{bertsi96a,bhabor99a,complearning99}.  
When the algorithm is linear,
then the  error equations  take 
the following linear  recursive form:
\begin{equation}
X_{t+1} =\left[ I-\alpha M_t \right] X_t
+ W_{t+1} , 
\elabel{RLSalpha}
\end{equation}
where $\bfmX=\{X_t\}$ is an error sequence,  
$\bfmM=\{M_t\}$ is a sequence of
$k \times k$ random matrices, $\bfmW=\{W_t\}$ 
is a ``disturbance'', and $I$ is the 
$k \times k$ identity matrix.  

An important example is the LMS
(least mean square) algorithm.
Consider the discrete linear time-varying model:
\begin{equation}
y(t)=\theta(t)^T \phi(t)+n(t), \hspace*{0.1in} t\geq 0
\label{model}
\end{equation}
where $y(t)$ and $n(t)$ are the sequences of
(scalar) observations and noise, respectively, 
and $\theta(t)=[\theta_1(t),\theta_2(t),\ldots ,\theta_k(t)]^T$ 
and $\phi(t)=[\phi_1(t),\ldots ,\phi_k(t)]^T$ 
denote the $k$-dimensional regression vector and 
time varying parameters, respectively.
The LMS algorithm is given by the recursion
\begin{equation}
\hat{\theta} (t+1)= \hat{\theta} (t)+\alpha \phi (t) e(t) ,
\end{equation}
where $e(t) \triangleq y(t)- \hat{\theta} (t)^T \phi (t)$, 
and the parameter $\alpha \in
(0,1]$ is the {\it step size}. Hence,
\begin{equation}
\tilde{\theta} (t+1) =(I-\alpha \phi (t) \phi
(t)^T )\tilde{\theta} (t) + \left[ \theta (t+1) - \theta (t) - \alpha \phi
    (t)n(t)  \right]   ,
\end{equation}
where $\tilde{\theta} (t) \triangleq \theta (t) - \hat
{\theta} (t)$.  This is of the form \eq RLSalpha/ with
$M_t=\phi (t) \phi(t)^T$, $W_{t+1}=\theta (t+1) - \theta (t) - \alpha
\phi (t)n(t)$, and $X_t=\tilde{\theta} (t)$.

On iterating \eq RLSalpha/  
we obtain the representation,
\be
\hspace*{-0.1in}
X_{t+1}
&=&\left( I-\alpha M_t \right) X_t + W_{t+1}	
	\nonumber\\
&=&\left( I-\alpha M_t \right) \left[ \left( I-\alpha M_{t-1} \right) 
  X _{t-1} + W_t \right] +W_{t+1}
	\label{e:EstErr}\\
&=&  \prod_{i=t}^0 (I-\alpha M_i)  X_0 + 
	\prod_{i=t}^1 (I-\alpha M_i)  W_1 + \dots  + 
	(I-\alpha M_t)W_t
	+W_{t+1}.
	\nonumber
\ee
From the last expression it is clear that the 
matrix products $\prod_{i=t}^s
(I-\alpha M_i)$  play an important role in the behavior of \eq RLSalpha/.

Properties of products of random matrices are of interest in a wide
range of fields. Application areas include numerical analysis 
\cite{ghaana99,vis00a},  statistical physics \cite{bou87a,cripalvul93a},  
recursive algorithms \cite{dabmas00,mou98a}, perturbation 
theory for dynamical systems \cite{arn98}, queueing theory \cite{mai97}, 
and even botany \cite{roe89a}.  Seminal results are contained in 
\cite{bel54a,ose68a,ose65a}.

A complementary and popular research area  concerns the eigenstructure
of \textit{large} random matrices (see e.g.\ \cite{tsehanste98a,hantse98a} for recent application to capacity of
communication channels). Although the results of the present paper do not 
address these issues, they provide justification for simplified  models
in communication theory, leading to bounds 
on the capacity for time-varying communication channels
\cite{medmeyhua01a}.

\notes{Jianyi to improve discussion on LMS and comm apps}

The relationship with dynamical systems theory is particularly relevant 
to the issues addressed here. Consider a {\em nonlinear}
dynamical system described by the equations,
\begin{equation}
X_{t+1} = X_t - \alpha f( X_t, \Phi_{t+1} ) + W_{t+1}\, ,
\elabel{Nonlin}
\end{equation}
where $\bfPhi=\{\Phi_t\}$ is an ergodic Markov process, evolving on a state
space $\state$, and $f\colon\Re^k \times\state \to\Re^k$ is smooth
and Lipschitz continuous. Although it is, 
of course, impossible to iterate a
nonlinear model of this general form, 
we can construct a random linear model 
to address many interesting issues.  
Viewing the initial condition $\gamma=X_0\in\Re^k$ 
as a continuous  variable, we write $X_t(\gamma)$ 
as the resulting state trajectory and consider 
the sensitivity matrix, 
\[
S_{t} = \frac{\partial}{\partial \gamma} X_t(\gamma),\qquad t\ge 0\, .
\]
From \eq Nonlin/ we have the linear recursion,
\begin{equation}
S_{t+1} =[I - \alpha M_{t+1}] S_t,
\elabel{sensitivity}
\end{equation}
where $M_{t+1} = \nabla_x f\, (X_t,\Phi_{t+1})$, $t\ge 0$.
If $\bfmS=\{S_t\}$ is suitably stable then the same is true for the nonlinear
model, and we find that trajectories couple to a steady state process
$\bfmX^*=\{X_t^*\}$: 
\[
\lim_{t\to\infty} \|X_t(\gamma) - X^*_t\| = 0\, . 
\]
These ideas are related to issues  developed in \Section{nonlin}.

The traditional analytic technique for addressing the 
stability of~\eq Nonlin/ or of~\eq RLSalpha/ is the 
\textit{ODE method} of \cite{lju77a}.  
For linear models, the basic idea
is that, for small values of $\alpha$, the behavior of
\eq RLSalpha/ should mimic that of the linear ODE,
\begin{equation}
\ddt\gamma_t=-\alpha
\barM \gamma_t+\barW\, ,
\elabel{ODE}
\end{equation}
where $\barM$ and $\barW$ are means of $M_t$ and
$W_t$, respectively.
To obtain a 
finer performance analysis one can instead compare  \eq RLSalpha/ 
to the linear diffusion model,
\begin{equation}
d \Gamma_t=-\alpha \barM \Gamma_t+dB_t,
\elabel{SDE}
\end{equation}
where $\bfmB=\{B_t\}$ is a Brownian Motion.

Under certain assumptions one may show that,
if the ODE \eq ODE/ is
stable, then the stochastic model \eq RLSalpha/ is stable 
in a statistical sense, and
comparisons with \eq SDE/ are possible under
still stronger assumptions (see e.g.\ \cite{ben99a,bormey00a,kusyin97a,kus84a}
for results concerning both linear and nonlinear recursions).

In  \cite{mou98a} an alternative point of view was proposed where
the stability verification problem
for \eq RLSalpha/ is cast in terms of the spectral radius of 
an associated discrete-time semigroup of linear operators.
This approach is based on the functional analytic setting of \cite{meytwe93e},
and analogous techniques are used in the treatment of multiplicative ergodic
theory and spectral theory in \cite{balmey98a,konmey01a,konmey01b}.
The main results of  \cite{mou98a} may be interpreted 
as a significant extension of the ODE method for linear 
recursions. 

\medskip

Our present results give a unified treatment of 
both the linear and nonlinear models treated in
\cite{mou98a} and \cite{bormey00a}, 
respectively.\footnote{Our results are given 
here with only brief proof outlines; 
a more detailed and complete account is 
in preparation.} Utilizing the operator-theoretic 
framework developed in \cite{konmey01a} also 
makes it possible to offer a transparent treatment,
and also significantly weaken the
assumptions used in earlier results.

We provide answers to the following questions:
\begin{description}
\item[(i)]
 For what range of $\alpha >0$ is the random linear system \eq RLSalpha/ $L_2$-stable,
 in the sense that
 $\Expect_x[\|X_t\|^2]$ is bounded in $t$?

\item[(ii)]
What does the averaged model \eq ODE/ tell us about the behavior of  
the original stochastic model?

\item[(iii)]
What is the impact of variability on \textit{performance} of recursive
algorithms?
\end{description}

\section{Linear Theory}
\slabel{prods}

In this section we develop stability theory and structural results for
the linear 
model \eq RLSalpha/ where $\alpha\ge 0$ is a fixed constant.

It is assumed that an underlying
Markov chain $\bfPhi$, with general state-space $\state$,
governs the statistics of \eq RLSalpha/ in the sense that
$\bfmM$ and $\bfmW$ 
are functions of the Markov chain:
\begin{equation}
 M_t = m(\Phi_t),\quad W_t = w(\Phi_t),\qquad t\ge 0\, .
\elabel{Assumptions}
\end{equation}
We assume that the entries of the $k\times k$-matrix valued function $m$
are bounded functions of $x\in\state$.  
Conditions on the vector-valued
function $w$ are given below.

We begin with some basic assumptions on $\bfPhi$,
required to construct a
linear operator with useful properties.

\subsection{Some spectral theory}

We assume throughout that the
Markov chain $\bfPhi$ is \textit{geometrically ergodic}  
or, equivalently, \textit{$V$-uniformly ergodic}.
This is equivalent to assuming the validity of the
following two conditions:
\begin{description}
\item\textit{Irreducibility \&\ aperiodicity:}  
There exists a $\sigma$-finite measure $\psi$ on the 
state space $\state$ such that, for any $x\in\state$
and any measurable $A\subset\state$ with $\psi(A)>0$,
\[
P^t(x,A)\eqdef \Prob \{ \Phi_t\in A \mid \Phi(0)=x\} > 0,\;\;\;
\hbox{for all sufficiently large $t>0$.}
\]
\item\textit{Geometric drift:}
There exists a  {\em Lyapunov function}  $V:\state\to[1,\infty)$, 
$ \gamma<1$, $b<\infty$,  $t_0\ge 1$, a `small set' $C$,  
and a `small measure' $\nu$, satisfying
\begin{equation}
\begin{array}{rcll}
PV\, (x) &\leq& \gamma V(x) + b\ind_C (x),\qquad & x\in\state
\\
P^{t_0} (x,\varble)  &\ge &\nu(\varble)       \qquad & x\in C
\end{array}
\elabel{Lyapunov} 
\end{equation} 
\end{description}
Under these assumptions it is known that 
$\bfPhi$ is ergodic and has a unique 
invariant probability measure $\pi$, 
to which it converges geometrically
fast, and without loss of generality
we can assume that $\pi(V^2)<\infty.$
For a detailed development of 
geometrically ergodic Markov processes 
see \cite{MT,meytwe93e,konmey01a}.   

We let $\LV$ denote the   set of measurable 
\textit{vector-valued} functions 
$g\colon\state\to \Co^k$ satisfying 
\[
\|g\|_V \eqdef \sup_{x\in\state}\frac{\|g(x)\|}{V(x)} < \infty\,,
\]
where $\|\cdot\|$ is the Euclidean norm on $\Co^k$, and $V\colon\state\to [1,\infty)$ 
is the Lyapunov function as above.  
For a linear operator $\clL\colon\LV\to\LV$ we define
the induced operator norm via 
\[
\lll\clL\lll_V\eqdef\sup \|\clL f\|_V/\|f\|_V
\]
where the supremum is over all non-zero $f\in \LV$. 
We say that $\clL$ is a bounded linear operator
if $\lll\clL\lll_V<\infty$, and its spectral radius is then given by
\begin{equation}
\xi\eqdef \lim_{t\to\infty} \bigl(\lll \clL^t \lll\bigr)^{1/t}
\elabel{sradius}
\end{equation}
The \textit{spectrum} $S({\clL })$
of the linear operator $\clL$ is 
$$S({\clL }) \eqdef \{z \in \mathbb{C}: \left(Iz-
  \clL  \right) ^{-1}  \hbox{does not exist as a bdd
  linear operator on } L_{\infty}^V \}.$$
If $\clL $ is a finite matrix, its spectrum is just the  
collection of all its eigenvalues. Generally, for the linear operators
considered in this paper, the dimension of $\clL $ and its spectrum will 
be infinite.

The family of linear operators  $\clL_{\alpha}\colon \LV\to\LV$,
$\alpha\in\RL$,
that will be used to analyze the recursion \eq RLSalpha/ are
defined by,
\begin{equation}
\begin{array}{rcl} 
{\clL_{\alpha}}f(x) 
&\eqdef& \Expect \left[ (I-\alpha m(\Phi_1))^\transpose f(\Phi_1) \mid
  \Phi(0)=x \right]
\\[.25cm]
&=& \Expect_x \left[(I-\alpha M_1)^\transpose f(\Phi_1) \right]\, ,
\end{array} 
\end{equation}
and we let $\xi_\alpha$ denote the spectral radius of $\clL_\alpha$.  

We assume throughout the paper 
that $m\colon\state\to \Re^{k\times k}$ 
is a bounded function. 
Under these conditions we obtain the following result as in \cite{mou98a}. 
\begin{theorem}  
\tlabel{mou98}  There exists $\alpha_0>0$ 
such that  for $\alpha\in (0,\alpha_0)$, $\xi_\alpha<\infty$, 
and $\xi_\alpha \in S(\clL_\alpha)$. 
\qed
\end{theorem}

To ensure that the recursion \eq RLSalpha/ is stable 
it is necessary that the spectral radius satisfy
$\xi_\alpha <1$. Under this condition it is obvious that the mean 
$\Expect[X_t]$ is 
uniformly bounded in $t$.    
The following result summarizes
additional conclusions obtained below.
\begin{theorem}
\tlabel{mainLinear}
Suppose that the eigenvalues of 
$\barM\eqdef \int m(x)\, \pi(dx)$ are all positive, 
and that $w^2\in\LV$, where the square is interpreted component-wise. 
Then, there exists a bounded open set $O\in \Re$ containing  
$ (0,\alpha_0)$, where $\alpha_0$ is given in \Theorem{mou98}, 
such that:
\begin{description}
\item[(i)]
For all $\alpha\in O$ we have $\xi_\alpha<1$ ,
and for any initial condition $\Phi_0=x\in\state$, 
$X_0=\gamma\in\Re^k$, 
\[
\Expect_x[\|X_t\|^2] \to \sigma^2_\alpha<\infty,
\;\;\mbox{geometrically fast, as $t\to\infty$.}
\]
\item[(ii)]
If $\bfPhi$ is stationary, then
for $\alpha\in O$ there exists a stationary process 
$\bfmX^\alpha$ such that for any initial condition 
$\Phi_0=x\in\state$, $X_0=\gamma\in\Re^k$,
\[
\Expect[\| X_t(\gamma) - X^\alpha_t\|^2]\to 0,
\;\;\mbox{geometrically fast, as $t\to\infty$.}
\]
\end{description}

\item[(iii)]
If $\alpha\not\in \barO$ and $\bfmW$ is i.i.d.\ with $\Sigma_W^2\neq 0$ 
then $ \Expect_x[\|X_t\|^2]$ is unbounded.
\end{theorem}

\begin{figure}[ht]
\begin{center} 
\Ebox{.5}{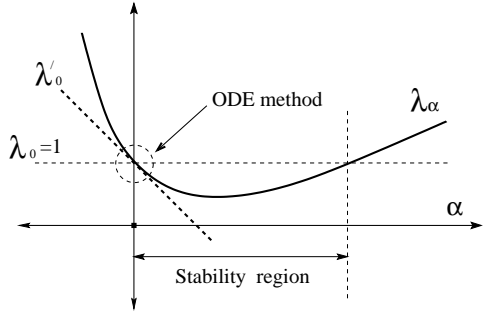}
\end{center}
\flabel{comparison}
\caption{The 
graph shows how $\la:=\xi_\alpha$ varies with $\alpha$. 
When $\alpha$ is close to $0$  \Th{operatorL} below implies
that the ODE \eq ODE/ determines stability of the algorithm since
it determines whether or not $\xi_\alpha<1$. A second-derivative
formula is also given in  \Th{operatorL}: If $\lambda''_0$ is large,
then the range of $\alpha$ for stability will
be correspondingly small.}
\end{figure}
\noindent
{\sc Proof Outline for \Theorem{mainLinear}}
Iterating the system equation \eq EstErr/ we may express
the expectation $\Expect_x \left[X_{t+1}^T X_{t+1} \right]$
as a sum of terms of the form,
\begin{equation}
\textstyle
\Expect_x \left[ W_j^T \Bigl(\prod_{i=t}^j (I-\alpha M_i)\Bigr)^T \Bigl(\prod_{i
=t}^k (I-\alpha M_i)\Bigr) W_k \right],\quad j,k=0,\dots , t\, .
\elabel{Qerr}
\end{equation}
For simplicity consider the case $j=k$.  Taking conditional expectations
at time $j$, one can then express the expectation \eq Qerr/ as
\[
\textstyle
\trace\Bigl(
\Expect_x\Bigl[ (\clQ_\alpha^{t-j} h\, (\Phi_j)) w(\Phi_j)w(\Phi_j)^T\Bigr]
\Bigr)
\]
where $\clQ_\alpha$ is defined in \eq clQ/, and $h\equiv I_{k\times k}$.
We define $O$ as the set of $\alpha$ such that the spectral radius of
this linear operator is strictly less than unity. Thus, for $\alpha\in O$
we have, for some   $\eta_\alpha<1$,
\[
\trace\Bigl(
 (\clQ_\alpha^{t-j} h\, (y)) w(y)w(y)^T
\Bigr)
=O( V(y)^2 e^{-\eta_\alpha (t-j)} ),\qquad \Phi_j=y\in\state.
\]
Similar reasoning may be applied for arbitrary $k,j$, and this shows
that $\Expect[\|X_t\|^2]$ is bounded  in $t\ge 0$
for any deterministic initial
conditions $\Phi_0=x\in\state$, $X_0=\gamma\in\Re^k$.
 
To construct the stationary process $X^\alpha$ we apply
backward coupling as developed in \cite{tho00}. 
Consider the system
starting at time $-n$, initialized at $\gamma=0$, and let $X^{\alpha,n}_t$,
$t\ge -n$, denote the resulting state trajectory.  We then have
for all $n,m\ge 1$,
\[
X^{\alpha,m}_t - X^{\alpha,n}_t
=  \Bigl(\prod_{i=t}^0 (I-\alpha M_i)\Bigr)[X^{\alpha,n}_0-X^{\alpha,n}_0 ],
\qquad t\ge 0\, ,
\]
which implies convergence in $L_2$ to a stationary process:  $X^\alpha_t\eqdef
\lim_{n\to\infty} X^{\alpha,n}_t$, $t\ge 0$.
We can then compare to the process initialized at  $t=0$,
\[
X^\alpha_t - X_t(\gamma) =  \Bigl(\prod_{i=t}^0 (I-\alpha M_i)\Bigr)[X^\alpha_0
- X_0(\gamma)],\qquad t\ge 0\, ,
\]
and the same reasoning as before gives (ii).
\qed

\subsection{Spectral decompositions}

Next we show that $\lambda_\alpha:=\xi_\alpha$ is 
in fact an eigenvalue of $\clL_\alpha$ for a range 
of $\alpha\sim 0$, and we use this fact to obtain 
a multiplicative ergodic theorem.
The maximal eigenvalue $\lambda_\alpha$ in 
\Theorem{evalue} is a generalization of the 
Perron-Frobenius eigenvalue;
c.f. \cite{sen80,konmey01a}.
\begin{theorem} 
\tlabel{evalue}
Suppose that the eigenvalues $\{\lambda_i(\barM)\}$ of $\barM $ are 
distinct.  Then,
\begin{description}
\item[(i)]
There exists $\epsy_0>0$ such that the linear operator $\clL_z$ has $k$ distinct
eigenvalues  $\{\lambda_{1,z},\dots,\lambda_{k,z}\}\subset S(\clL_z)$  for all $z\in B(\epsy_0) \eqdef \{ z\in\Co: |z-1|<\epsy_0 \}$, and 
$\lambda_{i,z}$ is an analytic
function of $z$ in this domain for each $i$.

\item[(ii)]
For $z\in B(\epsy_0)$ 
there are associated eigenfunctions  $\{h_{1,z},\dots,h_{k,z}\}\subset \LV$ 
and eigenmeasures  $\{\mu_{1,z},\dots,\mu_{k,z}\}\subset \clM_1^V$ satisfying 
\[
\clL_z h_{i,z} = \lambda_{i,z} h_{i,z} ,\qquad  \mu_{i,z} \clL_z   = \lambda_{i,z}\mu_{i,z}\, .
\]
Moreover, for each $i$, $x\in\state$, $A\in\bx$,  $\{h_{i,z}(x),\mu_{i,z}(A)\}$
are analytic functions on $B({\epsy_0})$.

\item[(iii)]
Suppose moreover that the eigenvalues $ \{ \lambda_i(\barM) \}$ are real.
Then we may take $\epsy_0>0$ sufficiently small so that 
$\{\lambda_{i,\alpha}, h_{i,\alpha},\mu_{i,\alpha}\}$ 
are real for $\alpha\in (0,\epsy_0)$. 
The maximal eigenvalue $\lambda_\alpha \eqdef \max_i \lambda_{i,\alpha}$ 
is equal to $\xi_\alpha$, and the corresponding eigenfunction and eigenmeasure
may be scaled so that the following limit holds:
\[
\lambda^{-t}_\alpha\clL^t_\alpha \to \ha\otimes\mu_\alpha,\qquad t\to\infty,
\]
where the convergence is in the $V$-norm.  

In fact,
there exists $\delta_0>0$ and $b_0<\infty$
such that for any $f\in\LV$ the following limit holds:
\[
\lambda_{\alpha}^{-t} \Expect_x \Bigl[ \Bigl(
  \prod_{i=1}^t (I-\alpha M_i)  \Bigr)^\transpose f(\Phi_t) \Bigr]
= h_{\alpha} (x) \mu_\alpha(f) + b_0 e^{-\delta_0 t} V(x)\, .
\]
\end{description}
\end{theorem}

\proof
The linear operator $\clL_0$ possesses a $k$-dimensional 
eigenspace corresponding to 
the eigenvalue $\lambda_0=1$.  This eigenspace is precisely the set of constant functions,
with a corresponding basis of eigenfunctions  given by  $\{e^i\}$,
where $e^i$ is the $i$th basis element in $\Re^k$.  The $k$-dimensional set of
vector-valued eigenmeasures $\{\pi^i\}$  given by 
$\pi^i = {e^i}^\transpose \pi$ spans the set
of all eigenmeasures with eigenvalue $\lambda_{0,i}=1$.

Consider the linear operator defined by 
\[
\Pi f\, (x) = (\pi(f_1),\dots,\pi(f_k))^\transpose
=  \bigl[ \sum e^i \otimes\pi^i \bigr] f
,\qquad f\in\LV.
\]  
It is obvious that  $\Pi\colon\LV\to\LV$ is a rank-$k$ linear operator, and 
for $\alpha =0$ we have from the $V$-uniform ergodic theorem of \cite{MT},
\[
\clL^t_0 - \Pi  = [\clL_0 - \Pi]^t \to 0,\qquad t\to\infty,
\] 
where the convergence is in norm, and hence takes place exponentially fast.
It follows that the spectral radius of 
$(\clL_0 - \Pi)$ is strictly less than unity.
By standard arguments it follows that, for some $\epsy_0>0$,
the spectral radius of $ \clL_z - \Pi$ is  also strictly less than unity.
The results then follow as in Theorem~3 of \cite{konmey01b}.
\qed

\medskip

Conditions under which the  bound  $\xi_\alpha <1$  
is satisfied are given in \Theorem{operatorL}, 
where we also 
provide formulae for the derivatives 
of $\lambda_\alpha$:
\begin{theorem} 
\tlabel{operatorL}
Suppose that the eigenvalues $ \{ \lambda_i(\barM) \}$ are real
and distinct.
Then, the maximal eigenvalue $\lambda_\alpha=\xi_\alpha$ satisfies,
\begin{description}
\item[(i)]
$ \frac{d}{d\alpha} \lambda_\alpha\Big|_{\alpha=0} =- \lmin(\barM)$.

\item[(ii)] 
The second derivative is given by,
\[
\frac{d^2}{d\alpha^2} \lambda_\alpha\Big|_{\alpha=0}
=  2\sum_{l=0}^{\infty}
v_0^\transpose\Expect_\pi[(M_0 - \barM)(M_{l+1} - \barM)] r_0 \, ,
\]
where $r_0$ is a right eigenvector of $\barM$ corresponding 
to $\lmin(\barM)$, and
$v_0$ is the left eigenvector, 
normalized so that $v_0^\transpose r_0=1$.

\item[(iii)]
Suppose that $m(x)=m^\transpose(x)$, $x\in\state$.  
Then we may take $v_0=r_0$
in (ii), and the second
derivative may be expressed, 
\[
\frac{d^2}{d\alpha^2} \lambda_\alpha\Big|_{\alpha=0} = \trace(\Gamma - \Sigma)
\]
where an $\Gamma$ is the  Central Limit Theorem
covariance for the stationary 
vector-valued stochastic process $F_k=[M_k-\barM] v_0$, 
and $\Sigma=\Expect_\pi[F_kF_k^\transpose]$ is its variance.  
\end{description}
\end{theorem}

\proof 
To prove (i), we differentiate the eigenfunction equation 
$\La \ha=\la\ha $ to obtain 
\begin{equation}
\La' \ha + \La \ha'=\la'\ha+\la\ha'.
\elabel{Lfirst}
\end{equation}
Setting $\alpha=0$ then gives a version of \textit{Poisson's equation},
\begin{equation}
{\clL}'_0 \ho +P\ho'=\lo'\ho+\ho',
\elabel{MatrixFish}
\end{equation}
where ${\clL}'_0 \ho= \Expect_x \left[-m(\Phi_1)^\transpose h_0(\Phi_1)\right]$. 
Since $\ho'\in\LV$ we may integrate both   sides with respect to the
invariant probability $\pi$ to obtain
\[
\Expect_\pi \left[-m(\Phi_1)^\transpose \right]h_0=
-\barM^\transpose h_0=\lo'\ho.
\] 
This shows that $\lo'$ is an eigenvalue of $-\barM$, and $\ho$ is an associated
eigenvector for $\barM^\transpose$.   It follows that 
$\lo'=-\lmin(\barM)$ by maximality of $\lambda_\alpha$.

We note that Poisson's equation \eq MatrixFish/ combined with
equation~(17.39) of \cite{MT} implies the formula,
\begin{equation} 
\ho'(x) 
= \Expect_{\pi}[\ho'(\Phi(0))]  
- \sum_{l=0}^{\infty}\Expect_{x}[(M_{l+1} - \barM)^\transpose] h_0 \, .
\end{equation}

To prove (ii) we consider the second-derivative formula,
\[
\clL_\alpha''\ha 
+ 2 \clL_\alpha' h'_\alpha
+ \clL_\alpha h''_\alpha  
=  
\lambda_\alpha h''_\alpha
+
2\lambda'_\alpha h'_\alpha
+
\lambda''_\alpha \ha.
\]
Evaluating these expressions at $\alpha=0$ 
and integrating with respect to $\pi$ then gives
the steady state expression,
\begin{equation}
\lambda''_0 h_0  
=
-2\Expect_\pi [ (M_1+ \lambda'_0) h'_0(\Phi_1)].
\elabel{pplambda}
\end{equation}
In deriving this identity we have used the expressions,
\[
\clL_0'f\, (x) = \Expect_x[M_1 f(\Phi_1)],\quad
\clL_0''f\, (x) = 0,\qquad f\in\LV,\  x\in\state.
\]
This combined with \eq pplambda/ gives the desired formula
since we may take $v_0=h_0$ in (ii).

To prove (iii) we simply note that in the symmetric case the formula
in (ii) becomes,
\[
\lambda''_0 = \sum_{k\neq 0} \Expect_\pi[\|F_k\|^2] 
        = \trace(\Gamma - \Sigma)\, .
\]

\subsection{Second-order statistics}

In order to understand the second-order statistics of
$\bfmX$ it is convenient to introduce 
another linear operator $\clQ_\alpha$ as follows,
\begin{equation}
\begin{array}{rcl} 
{\clQ_{\alpha}}f(x) &=& \Expect \left[ (I-\alpha m(\Phi_1))^\transpose f(\Phi_1) (I-\alpha m(\Phi_1))|
  \Phi(0)=x \right]
\\[.25cm]
&=& \Expect_x \left[(I-\alpha M_1)^\transpose f(\Phi_1) (I-\alpha M_1) \right],
\end{array} 
\elabel{clQ}
\end{equation}
where the domain of $\clQ_{\alpha}$
is the collection of matrix-valued functions 
$f\colon\state\to \Co^{k\times k}$.
When considering $\clQ_\alpha$ we redefine $L_{\infty}^V$ accordingly. 
It is clear that  $\clQ_\alpha:L_{\infty}^V \rightarrow L_{\infty}^V $ 
is a bounded linear operator
under the geometric drift condition
and the boundedness assumption on $m$.

Let $\xi^Q_z$ denote the spectral radius of
$\clQ_\alpha$. We can again argue that $\xi^Q_z$ 
is smooth in a neighborhood of the origin, 
and the following follows as in \Theorem{operatorL}:
\begin{theorem} 
\tlabel{operatorQ} 
Assume that the eigenvalues of $\barM$ are real and distinct.
Then there exists $\epsy_0>0$ such that for each $z\in B({\epsy_0})$
there exists an eigenvalue $\eta_z\in\Co$ for $\clQ_z$
satisfying $|\eta_z| = \xi^Q_z$, and $\eta_\alpha$ is real
for real $\alpha\in (0,\epsy_0)$. The eigenvalue $\eta_z$ is smooth
on $B({\epsy_0})$ and satisfies,
\[
\eta_0'({\clQ}) = - 2\lmin(\barM).
\]
\end{theorem}

\proof
This is again based on differentiation of the eigenfunction equation
given by $\Qa \ha=\eta_\alpha \ha$, where $\eta_\alpha$ and $\ha$ are the
eigenvalue and matrix-valued eigenfunction, respectively.
Taking derivatives on both sides gives
\begin{equation}
\Qa' \ha + \Qa \ha'=\eta'_\alpha \ha+\eta_\alpha \ha'
\end{equation}
where ${\clQ}'_0 \ho= \Expect_x \left[-m(\Phi_1)^\transpose
h_0(\Phi_1)-h_0(\Phi_1)m(\Phi_1)\right]$. As before, we then
obtain the steady-state expression,
\begin{equation}
\Expect_\pi \left[-m(\Phi_1)^\transpose h_0-h_0
  m(\Phi_1)\right]=-\barM^\transpose h_0-h_0 \barM=\eta'_0 \ho.
\end{equation}
And, as before, we may conclude that $\eta_0'=2\lambda_0'=- 2\lmin(\barM)$.
\qed

\subsection{An illustrative example}

Consider the discrete-time, linear time-varying model
\begin{equation}
y_t=\theta_t^\transpose \phi_t+N_t, \qquad t\geq 0\, ,
\elabel{model}
\end{equation}
where $\bfmy$ is a sequence of
scalar observations, $\bfmN=\{N_t\}$ is a noise process, 
$\{\theta_t\}$ is the sequence of $k$-dimensional 
regression vectors, and $\{\phi_t\}$ are 
$k$-dimensional time-varying parameters.
In this section we illustrate the results above using the LMS 
(least mean square) parameter estimation algorithm,
\[
\hatheta _{t+1}= \hatheta _t+\alpha \phi _t e_t ,
\]
where $\bfme$ is the error sequence,
$e_t \eqdef y_t- \hatheta _t^\transpose \phi _t$, $t\ge 0$. 

As in the Introduction,
writing $\tiltheta_t = \theta_t-\hatheta_t$ we obtain
\[
\tilde{\theta} _{t+1} =(I-\alpha \phi _t \phi
_t^\transpose )\tilde{\theta} _t + \left[ \theta _{t+1} - \theta _t - \alpha \phi 
    _t N_t  \right]   \, .
\]
This is of the form \eq RLSalpha/ with
$X_t=\tilde{\theta} _t$, $M_t=\phi _t \phi_t^\transpose$ and 
$W_{t+1}=\theta _{t+1} - \theta _t -\alpha \phi _tN_t$.
 
For the sake of simplicity 
and to facilitate explicit numerical calculations,
we consider the following special case:
We assume that $\bfphi$ is of the form  
$\phi_t=(s_t,s_{t-1})^\transpose$, where  
the sequence $\bfms$ is Bernoulli ($s_t=\pm 1$
with equal probability) 
and take $\bfmN$ to be an i.i.d.
noise sequence.

In analyzing the random linear system 
we may ignore the noise $\bfmN$ and take $\bfPhi = \bfphi$.
This  is clearly geometrically ergodic 
since it is an ergodic, finite state space Markov chain,
with four possible states. In fact, $\bfPhi$ is
geometrically
ergodic with Lyapunov function
$V\equiv 1$. Viewing $h \in L_{\infty}^{V}$ 
as a vector in $\RL^8$, 
the eigenfunction equation for
$\clL_\alpha$ becomes 
\begin{equation}
\La \ha=\frac{1}{2}
\left[ \begin{array}{cccc}
A_1 & A_0 & A_2 & A_0 \\
A_1 & A_0 & A_2 & A_0 \\
A_0 & A_2 & A_0 & A_1 \\
A_0 & A_2 & A_0 & A_1
\end{array} \right] \ha
= \la\ha 
\end{equation}
where 
$A_0=\left[ \begin{array}{cc}
0 & 0 \\
0 & 0 
\end{array} \right]$,
$A_1=\left[ \begin{array}{cc}
1-\alpha & -\alpha \\
-\alpha  & 1-\alpha 
\end{array} \right]$,
$A_2=\left[ \begin{array}{cc}
1-\alpha & \alpha \\
\alpha & 1-\alpha
\end{array} \right]$.

\begin{figure}[ht]
\begin{center}
 \parbox{.48\hsize}{\epsfxsize=\hsize \epsfbox{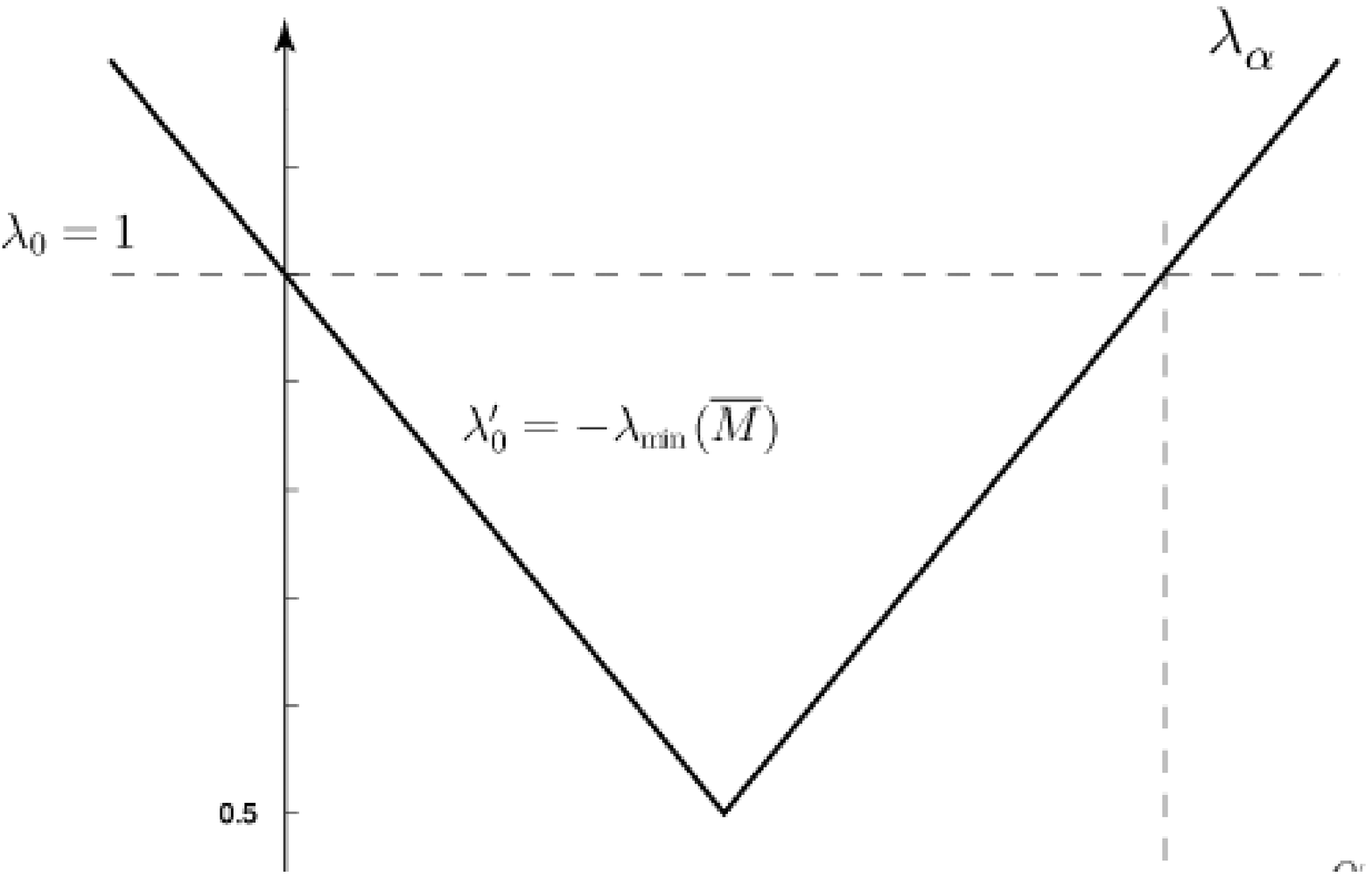}}
 \parbox{.48\hsize}{\epsfxsize=\hsize \epsfbox{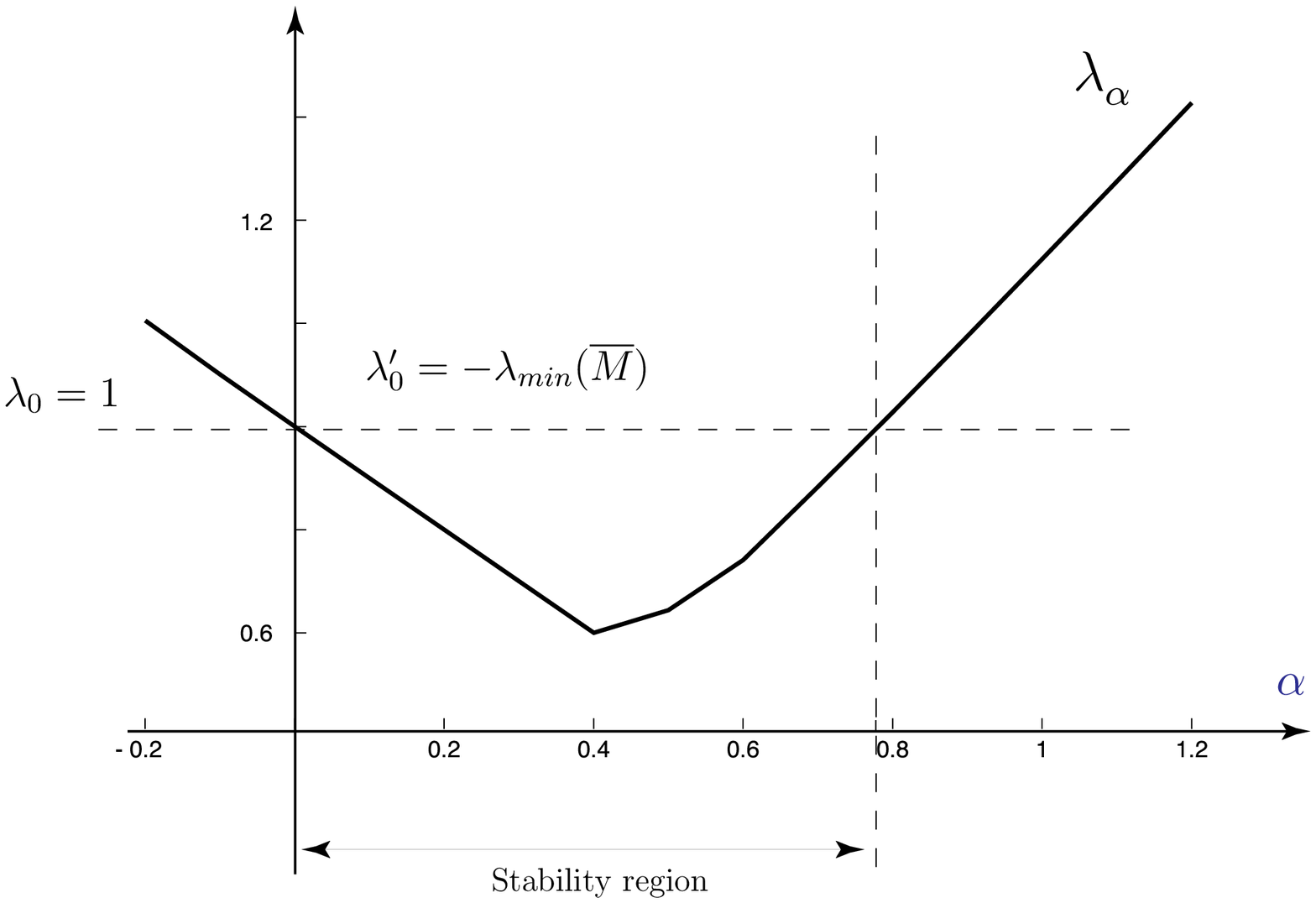}}
\end{center}
\caption{The figure on the left shows the Perron-Frobenius eigenvalue
$\lambda_\alpha= \xi_\alpha$ for the LMS model with 
$\phi_t=(s_t,s_{t-1})^\transpose$. The figure on the
right shows the case where $\phi_t=(s_t,s_{t-1},s_{t-2})^\transpose$.
In both cases, the sequence $\bfms$ is i.i.d.\ Bernoulli.}
\flabel{LK2} 
\Ebox{.65}{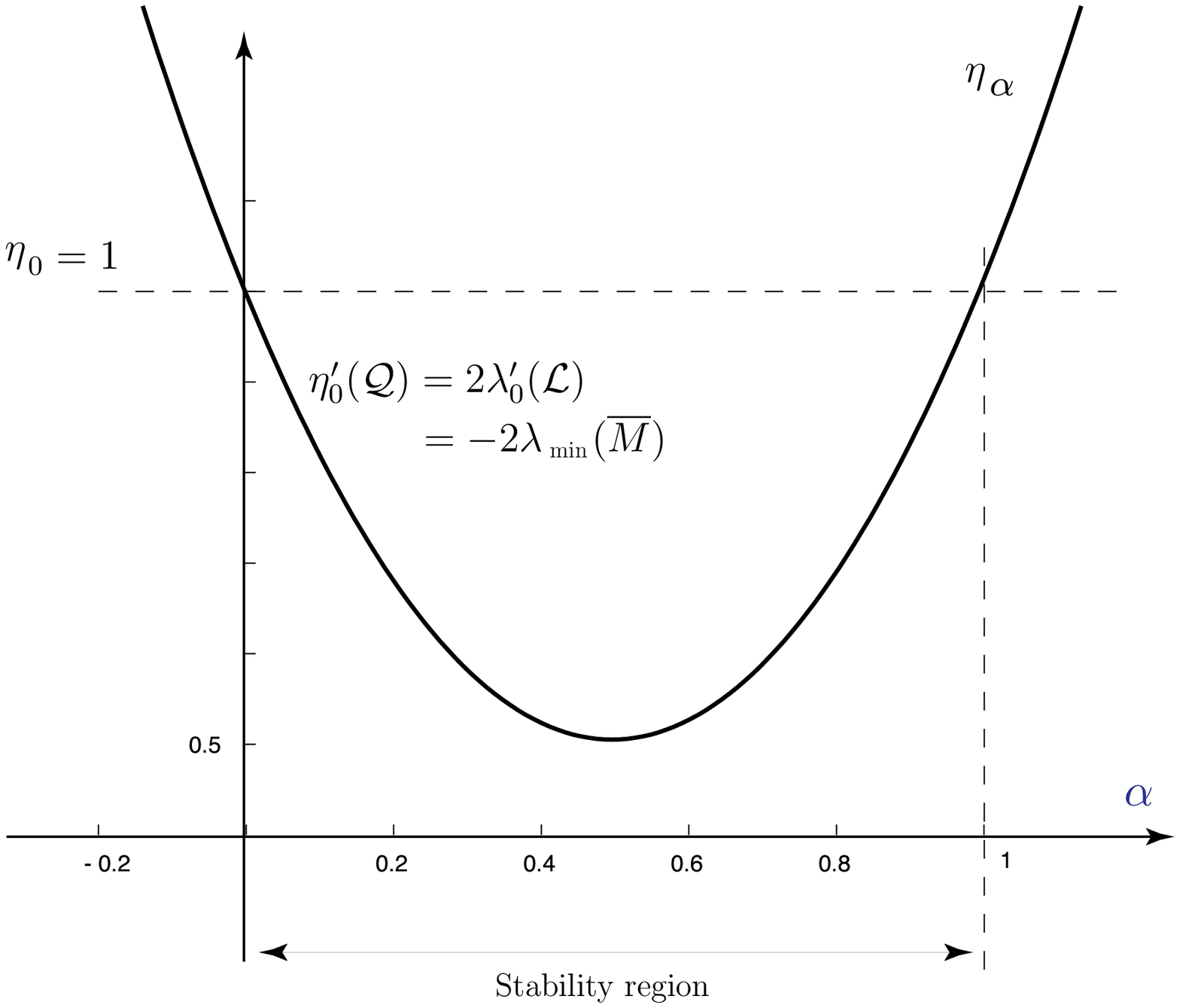}
\caption{The maximal eigenvalues $\eta_\alpha = \xi^Q_\alpha$ 
are piecewise quadratic in $\alpha$ in the case 
where $\phi_t=(s_t,s_{t-1})^\transpose$ with $\bfms$ as above.}
\flabel{Q_piecewise_linear}
\end{figure}

In this case, we have the following local behavior:
\begin{theorem}
\tlabel{exLocal}
In a neighbor of $0$, the spectral radii of  $\clL_\alpha$,  $\clQ_\alpha$
satisfy
\[
\begin{array}{rclrcl}
 \frac{d}{d\alpha} \xi_\alpha\Big|_{\alpha=0} &= &- \lmin(\barM) = -1;  
	\qquad &
  	\frac{d}{d\alpha} \xi^Q_\alpha\Big|_{\alpha=0} &= &- 2\lmin(\barM) = -2
\\[.25cm]
 \frac{d^n}{d\alpha^n} \xi_\alpha\Big|_{\alpha=0} &= & 0, n \geq 2;  
	\qquad &
  	\frac{d^n}{d\alpha^n} \xi^Q_\alpha\Big|_{\alpha=0} &= & 0, n \geq 3.
\end{array}
\]
So $\lambda_\alpha$ and $\eta_\alpha$ are linear 
and quadratic around $0$, respectively.
\end{theorem}

\proof
This follows from differentiating the respective eigenfunction
equations. Here we only show the proof for operator $\clQ$, the proof for
operator $\clL$ is similar.
 
Taking derivatives on both sides of the eigenfunction equation for $\clQ_\alpha$ gives,
\begin{equation}
\Qa' \ha + \Qa \ha'=\eta_\alpha'\ha+\eta_\alpha\ha'
\elabel{Qfirst}
\end{equation}
Setting $\alpha=0$ gives a version of \textit{Poisson's equation},
\begin{equation}
{\clQ}'_0 \ho + {\clQ} \ho'=\eta_0'\ho+\eta_0\ho'
\elabel{Qfirst_0}
\end{equation}

Using the identities of $h_0$ and  
${\clQ}'_0 h_{0} = \Expect_x[-M_1^\transpose h_0 -h_0 M_1]$, we obtain the
steady state expression
\begin{equation}
\barM^\transpose h_0+h_0\barM=-\eta'_0 h_0.
\elabel{Qfirst_steady}
\end{equation}
Since $\barM=I$, we have $\eta'_0=-2$.
Now, taking the 2nd derivatives on both sides of
\eq Qfirst/ gives,
\begin{equation}
{\clQ}''_\alpha \ha + 2{\clQ}'_\alpha
\ha'+{\clQ}_\alpha\ha''=\eta''_\alpha\ha+2\eta'_\alpha\ha'+\eta_\alpha\ha''.
\elabel{Qpprime}
\end{equation}
Letting $\alpha=0$ and considering the steady state, we obtain
\begin{equation}
2\barM^\transpose\ho\barM-2\Expect_{\pi}[M_1^\transpose\ho'+\ho'M_1]=\eta''_0\ho+2\eta'_0
\Expect_{\pi}[\ho'].
\elabel{Qsecond_0}
\end{equation}
Poisson's equation \eq Qfirst_0/ combined with equation \eq
Qfirst_steady/ and equation~(17.39) of \cite{MT} implies the formula,
\begin{equation}
\begin{array}{rcl}
\ho'(x) &=& \Expect_{\pi}(\ho')+\sum_{l=0}^{\infty}\Expect_{x}[-M_{l+1}^\transpose\ho-\ho 
M_{l+1}-\eta'_0\ho]
\\
&=& \Expect_{\pi}(\ho')+\sum_{l=0}^{\infty}\Expect_{x}[(\barM-M_{l+1})^\transpose\ho+\ho 
(\barM-M_{l+1})].
\end{array}
\end{equation}
So, from $\barM=I$, $\eta_0'=-2$ and \eq Qsecond_0/ we  have $\eta_0''=2$.
In order to show $\eta_\alpha$ is quadratic near zero, we take the 3rd
derivative on both sides of \eq Qpprime/ and consider the steady state 
at $\alpha=0$,
\begin{equation}
{\clQ}'''_0 \ho + 3{\clQ}''_0
\ho'+3{\clQ}'_0\ho''+{\clQ}_0\ho'''=\eta'''_0\ho+3\eta''_0\ho'+3\eta'_0\ho''+\eta_0\ho'''.
\elabel{Qpprime2}
\end{equation}
With equation~(17.39) of \cite{MT} and $\eta_0'=-2$ and $\eta_0''=2$,
we can show
$\eta_0'''=0$ and $\eta_0^{(n)}=0$ for $n>3$, hence $\eta_\alpha$ is
quadratic around $0$. 
\qed 

\section{Nonlinear models}
\slabel{nonlin}

We now turn to the nonlinear model shown in \eq Nonlin/.   
We take the special form,
\begin{equation}
X_{t+1} = X_t - \alpha [f( X_t, \Phi_{t+1} ) +  W_{t+1}]\, ,
\elabel{Nonlin2}
\end{equation}
We continue to assume that $\bfPhi$ is geometrically ergodic, and
that $W_t=w(\Phi_t)$, $t\ge 0$, with $w^2\in\LV$.  The associated 
ODE is given by
\begin{equation}
\ddt \gamma_t = \barf(\gamma_t),
\elabel{nonlinODE}
\end{equation}
where $\barf(\gamma) = \int f(\gamma,x)\, \pi(x)$, $\gamma\in\Re^k$.

We assume that 
$\barW = \Expect_\pi[W_1]=0$, and the following conditions are imposed 
on $f$:
\begin{description}
\item[(N1)] 
The function $f$ is Lipschitz, and there exists a function 
$\barf_\infty : \Re ^k  \to \Re ^k$ such that 
\[
\lim_{r \to \infty } r^{-1}\barf (r\gamma) = \barf_\infty (\gamma) , 
                \qquad \gamma \in \Re^k .
\]
Furthermore, the origin in $\Re ^k$ is an asymptotically
stable equilibrium point for the ODE,
\begin{equation}
\ddt \gamma^\infty_t = \barf_\infty (\gamma^\infty_t).
\elabel{ODEinfty} 
\end{equation}

\item[(N2)] 
There exists $b_f<\infty$ such that
$\displaystyle
\sup_{\gamma\in\Re^k} \|f(\gamma,x)-\barf(\gamma)\|^2 
        \le b_f V(x),\; x\in\state$.

\item[(N3)]
There exists a unique stationary point $x^*$ for the ODE \eq nonlinODE/
that is a globally asymptotically stable equilibrium.
\end{description}

Define the absolute error by
\begin{equation}
\epsy_t\eqdef\| X_t-x^*\|,    \qquad t\ge 0.
\elabel{error}
\end{equation}
The following result is an extension of Theorem~1 of \cite{bormey00a}
to Markov models:
\begin{theorem}
\tlabel{BSbdds}
Assume that (N1)--(N3) hold. Then there exists $\epsy_0>0$ such that 
for any $0<\alpha<\epsy_0$:
\begin{description}
\item[(i)]
For any $\delta > 0$, there exists $b_1=b_1(\delta)<\infty$ such that 
\[
\limsup_{n \to  \infty }  \Prob(\epsy_n\ge \delta ) \leq b_1 \alpha.
\]

\item[(ii)]
If the origin is a
globally exponentially asymptotically stable equilibrium for 
the ODE \eq nonlinODE/,
then there exists $b_2<\infty$ such that for every initial
condition $\Phi_0=x\in\state$, $X_0=\gamma\in\Re^k$,
\[
\limsup_{n\to\infty} \Expect[ \epsy_n^2 ]  \le b_2 \alpha.
\]
\end{description}
\end{theorem}

\noindent
{\sc Proof Outline for \Theorem{BSbdds}} The continuous-time process  
$\{x^o_t : t\ge 0\}$ is defined to be
the interpolated version of
$\bfmX$ given as follows: Let $T_j = j\alpha$, $j\ge 0$, 
and define $x^o(T_j) = \alpha X_j$, with $\bfmx^o$  defined by
linear  interpolation on the remainder of $[T(j), T(j+1)]$ to
form a piecewise linear function.
Using geometric ergodicity we can bound the error between $\bfmx^o$ and
solutions to the ODE \eq nonlinODE/ as in \cite{bormey00a}, and we may
conclude that the joint process $(\bfmX,\bfPhi)$ is 
geometrically ergodic with Lyapunov function
$V_2(\gamma,x) = \|\gamma\|^2 + V(x)$.
\qed

\medskip

We conclude with an extension of \Theorem{mainLinear}
describing the behavior of the sensitivity 
process $\bfmS$.   
\begin{theorem}
\tlabel{BSsens}
Assume that (N1)--(N3) hold, and that the eigenvalues of
the matrix $\barM$ have strictly positive real part, where
\[
\barM \eqdef  \nabla \barf\, (x^*)\,.
\]
Then there exists $\epsy_1>0$ such that 
for any $0<\alpha<\epsy_1$, the conclusions of 
\Theorem{BSbdds}~(ii) hold, and, in addition:
\begin{description}
\item[(i)]
The spectral radius $\xi_\alpha$ 
of the random linear system \eq sensitivity/ describing the evolution
of the sensitivity process is strictly less than
one. 

\item[(ii)]
There exists  a stationary process $\bfmX^\alpha$ such that for any initial
condition $\Phi_0=x\in\state$, $X_0=\gamma\in\Re^k$,
\[
\Expect[\| X_t - X^\alpha_t\|^2]\to 0,\qquad t\to\infty\, .
\]
\end{description}
\end{theorem}

\end{document}